\documentclass{amsart}
\usepackage{amsmath,amsthm,amssymb}

\makeatletter
    
    \@addtoreset{equation}{section}
  \makeatother

\newtheorem{definition}{Definition}[section]
\newtheorem{proposition}[definition]{Proposition}
\newtheorem{theorem}[definition]{Theorem}

\newtheorem{remark}{Remark}[section]

\pagebreak

\title[Lifespan of solutions to the semilinear damped wave equation]{
A note on the lifespan of solutions to the semilinear damped wave equation
}
\author[M. Ikeda and Y. Wakasugi]{Masahiro IKEDA and Yuta WAKASUGI}
\email{m-ikeda@cr.math.sci.osaka-u.ac.jp}
\email{y-wakasugi@cr.math.sci.osaka-u.ac.jp}
\address{Department of Mathematics, Graduate School of Science, Osaka University,
Toyonaka, Osaka, 560-0043, Japan}
\begin{document}
\begin{abstract}
This paper concerns estimates of the lifespan of solutions to
the semilinear damped wave equation
$\square u+\Phi(t,x)u_t=|u|^p$
in $(t,x)\in [0,\infty)\times\mathbf{R}^n$,
where the coefficient of the damping term is
$\Phi(t,x)=\langle x\rangle^{-\alpha}(1+t)^{-\beta}$
with $\alpha\in [0,1),\ \beta\in (-1,1)$
and
$\alpha\beta=0$.
Our novelty is to prove an upper bound of the lifespan of solutions
in subcritical cases $1<p<2/(n-\alpha)$.

\end{abstract}

\maketitle
\section{Introduction}
We consider the semilinear damped wave equation
\begin{equation}
\label{eq11}
	u_{tt}-\Delta u+\Phi(t,x)u_t=|u|^p,\quad (t,x)\in [0,\infty)\times \mathbf{R}^n,
\end{equation}
with the initial condition
\begin{equation}
\label{eqini}
	(u, u_t)(0,x)=\varepsilon (u_0, u_1)(x),\quad x\in\mathbf{R}^n,
\end{equation}
where $u=u(t,x)$ is a real-valued unknown function of $(t,x)$,
$1<p$,
$(u_0, u_1)\in H^1(\mathbf{R}^n)\times L^2(\mathbf{R}^n)$
and $\varepsilon$ is a positive small parameter.
The coefficient of the damping term is given by
$$
	\Phi(t,x)=\langle x\rangle^{-\alpha}(1+t)^{-\beta}
$$
with
$\alpha\in [0,1),\ \beta\in (-1,1)$
and
$\alpha\beta=0$.
Here $\langle x\rangle$ denotes $\sqrt{1+|x|^2}$.

Our aim is to obtain an upper bound of the lifespan of solutions to
(\ref{eq11}).

We recall some previous results for (\ref{eq11}).
There are many results
about global existence of solutions for (\ref{eq11})
and many authors have tried to determine the critical exponent
(see \cite{HKI, HO, IMN, ITY, LiZ, LNZ, Na, N1, O, TY1, W2, Z, Zhou05} and the references therein).
Here ``critical'' means that
if $p_c<p$, all small data solutions of (\ref{eq11}) are global;
if $1<p\le p_c$, the local solution cannot be extended globally even for small data.

In the constant coefficient case $\alpha=\beta=0$,
Todorova and Yordanov \cite{TY1} and Zhang \cite{Z} determined
the critical exponent of (\ref{eq11}) with compactly supported data as
$$
	p_c=p_F=1+\frac{2}{n}.
$$
This is also the critical exponent of the corresponding heat equation
$-\Delta v+v_t=|v|^p$
and called the Fujita exponent (see \cite{F}).

On the other hand, there are few results about upper estimates of the lifespan for (\ref{eq11}).
When $n=1, 2$, Li and Zhou \cite{LiZ}
obtained the sharp upper bound:
\begin{equation}
\label{eq14}
	T_{\varepsilon}\le \left\{\begin{array}{ll}
		\exp(C\varepsilon^{-2/n}),&\text{if}\ p=1+2/n,\\
		C\varepsilon^{-1/\kappa},&\text{if}\ 1<p<1+2/n,
	\end{array}\right.
\end{equation}
where $C=C(n,p,u_0,u_1)>0$
and $\kappa=1/(p-1)-n/2$
for the data $u_0, u_1\in C_0^{\infty}(\mathbf{R}^n)$ satisfying $\int (u_0+u_1)dx>0$.
Nishihara \cite{N2} extended this result to $n=3$
by using the explicit formula of the solution to the linear part of (\ref{eq11})
with initial data $(0 ,u_1)$:
$$
	u(t,x)=e^{-t/2}W(t)u_1+J_0(t)u_1.
$$
Here $W(t)u_1$ is the solution of the wave equation $\square u=0$ with initial data $(0, u_1)$
and $J_0(t)u_1$ behaves like a solution of the heat equation $-\Delta v+v_t=0$.
However,
both the methods of \cite{LiZ} and \cite{N2} do not work in higher dimensional cases $n\ge 4$,
because they used the positivity of $W(t)$,
which is valid only in the case $n\le 3$.
In this paper we shall extend both of the results to $n\ge 4$
in subcritical cases $1<p<1+2/n$.

Next, we recall some results of variable coefficient in cases $\alpha\neq 0$ or $\beta\neq 0$.
There are many results on asymptotic behavior of solutions
in connection with the diffusion phenomenon,
Here the diffusion phenomenon means that
solution of the damped wave equation behaves like
a solution for the corresponding heat equation as $t\rightarrow +\infty$.
For more details about the diffusion phenomenon, see, for example \cite{TY2, Wirth1, Wirth2}.

For the case $\alpha\in [0,1), \beta=0$,
Ikehata, Todorova and Yordanov \cite{ITY}
determined the critical exponent for (\ref{eq11}) as $p_c=1+2/(n-\alpha)$,
which also agrees
with that of the corresponding heat equation
$-\Delta v+\langle x\rangle^{-\alpha}v_t=|v|^p$.
Here we emphasize that in this case there are no results about upper estimates for the lifespan.
It will be given in this paper.

Next,
for the case $\beta\in (-1,1), \alpha=0$,
Nishihara \cite{N3} and Lin, Nishihara and Zhai \cite{LNZ}
proved $p_c=1+2/n$,
which is also same as that of the heat equation
$-\Delta v+(1+t)^{-\beta}v_t=|v|^p$.
On the other hand,
upper estimates of the lifespan have not been well studied.
Recently, Nishihara \cite{N3} obtained
a similar result of \cite{LiZ, N2}:
let $n\ge 1, \beta\ge 0$ and $(u_0,u_1)$ satisfy
$\int_{\mathbf{R}^n}u_i(x)dx\ge 0\ (i=0,1),\ \int_{\mathbf{R}^n}(u_0+u_1)(x)dx>0$.
Then there exists a constant $C>0$ such that
$$
	T_{\varepsilon}\le\left\{\begin{array}{ll}
		e^{C\varepsilon^{-(1+\beta)/n}},&\text{if}\ p=1+(1+\beta)/n,\\
		C\varepsilon^{-1/\hat{\kappa}},&\text{if}\ 1+2\beta /n\le p<1+(1+\beta)/n,
	\end{array}\right.
$$
where $\hat{\kappa}=(1+\beta)/(p-1)-n$.
We note that
the rate $\hat{\kappa}$ is not optimal,
because it is not same as that of the corresponding heat equation.
Moreover, there are no results
for $1+(1+\beta)/n<p\le 1+2/n$.
We note that the proof by Todorova and Yordanov \cite{TY1}
also gives the same upper bound in the case $\beta=0, 1<p<1+1/n$.
In this paper we will improve the above result for all $1<p<1+2/n$
and give the sharp upper estimate.

Finally, we mention that our method is not applicable to $\alpha\beta\neq 0$.
On the other hand,
the second author \cite{W2} proved a small data global existence result
for (\ref{eq11}) with $\alpha,\beta\ge 0, \ \alpha+\beta\le 1$,
when $p>1+2/(n-\alpha)$.
This also agrees with the critical exponent of
the corresponding heat equation $-\Delta v+\langle x\rangle^{-\alpha}(1+t)^{-\beta}v_t=|v|^p$.
Therefore, it is expected that when $1<p\le 1+2/(n-\alpha)$,
there is a blow-up solution for (1.1) in this case.

\section{Main Result}

First,
we define the solution of (\ref{eq11}).
We say that $u\in X(T):=C([0,T);H^1(\mathbf{R}^n))\cap C^1([0,T);L^2(\mathbf{R}^n))$
is a solution of (\ref{eq11}) with initial data (\ref{eqini}) on the interval $[0,T)$ if
the identity
\begin{align}
\label{sol}
	\lefteqn{\int_{[0,T)\times \mathbf{R}^n}
		u(t,x)(\partial_t^2\psi(t,x)-\Delta\psi(t,x)-\partial_t(\Phi(t,x)\psi(t,x)))dxdt}\notag\\
	&=\varepsilon\int_{\mathbf{R}^n}\left\{(\Phi(0,x)u_0(x)+u_1(x))\psi(0,x)-u_0(x)\partial_t\psi(0,x)\right\}dx\\
	&\quad+\int_{[0,T)\times\mathbf{R}^n}|u(t,x)|^p\psi(t,x)dxdt\notag
\end{align}
holds for any $\psi\in C_0^{\infty}([0,T)\times \mathbf{R}^n)$.
We also define the lifespan for the local solution of (\ref{eq11})-(\ref{eqini}) by
\begin{align*}
	T_{\varepsilon}:=\sup\{\, &T\in (0,\infty] ; 
		\text{there exists a unique solution $u\in X(T)$ of (\ref{eq11})-(\ref{eqini})}\,\}.
\end{align*}
We first describe the local existence result.
\begin{proposition}
Let $\alpha\ge 0, \beta\in\mathbf{R}, 1<p\le n/(n-2)\ (n\ge 3)$, $1<p<\infty \ (n=1,2)$,
$\varepsilon>0$ and $(u_0, u_1)\in H^1(\mathbf{R}^n)\times L^2(\mathbf{R}^n)$.
Then $T_{\varepsilon}>0$, that is,
there exists a unique solution
$u\in X(T_{\varepsilon})$
to $(\ref{eq11})$-$(\ref{eqini})$.
Moreover, if $T_{\varepsilon}<+\infty$, then it follows that
$$
	\lim_{t\to T_{\varepsilon}-0}\|(u, u_t)(t,\cdot)\|_{H^1\times L^2}=+\infty.
$$
\end{proposition}
For the proof, see, for example \cite{IT}.
Next, we give an alomost optimal lower estimate of $T_{\varepsilon}$.
\begin{proposition}
Let $(u_0, u_1)\in H^1(\mathbf{R}^n)\times L^2(\mathbf{R}^n)$ be compactly supported
and $\delta$ any positive number.
We assume that $\alpha\in [0,1), \beta\in (-1,1), \alpha\beta\ge 0$ and $\alpha+\beta<1$.
Then there exists a constant $C=C(\delta, n,p,\alpha, \beta, u_0,u_1)>0$ such that
for any $\varepsilon>0$, we have
$$
	C\varepsilon^{-1/\kappa+\delta}\le T_{\varepsilon},
$$
where
$$
	\kappa=\frac{2(1+\beta)}{2-\alpha}\left(\frac{1}{p-1}-\frac{n-\alpha}{2}\right).
$$
\end{proposition}
The proof of this proposition follows from the a priori estimate for the energy of solutions.
For the proof, see \cite{LNZ, ITY, W2}.
We note that the above proposition is valid even for the case $\alpha\beta\neq 0$.

Next, we state our main result,
which gives an upper bound of $T_{\varepsilon}$.
\begin{theorem}
Let $\alpha\in[0,1), \beta\in (-1,1), \alpha\beta=0$ and let $1<p<1+2/(n-\alpha)$.
We assume that the initial data $(u_0, u_1)\in H^1(\mathbf{R}^n)\times L^2(\mathbf{R}^n)$ satisfy
\begin{equation}
\label{eq16}
	\langle x\rangle^{-\alpha}Bu_0+u_1\in L^1(\mathbf{R}^n)\ \ \text{\and}\quad
	\int_{\mathbf{R}^n}(\langle x\rangle^{-\alpha}Bu_0(x)+u_1(x))dx>0,
\end{equation}
where
$$
	B=\left(\int_0^{\infty}e^{-\int_0^t(1+s)^{-\beta}ds}dt\right)^{-1}.
$$
Then there exists $C>0$ depending only on $n, p, \alpha, \beta$ and $(u_0, u_1)$ such that
$T_{\varepsilon}$ is estimated as
$$
	T_{\varepsilon}\le C\begin{cases}
		\varepsilon^{-1/\kappa}&\text{if}\quad 1+\alpha/(n-\alpha)<p<1+2/(n-\alpha),\\
		\varepsilon^{-(p-1)}(\log(\varepsilon^{-1}))^{p-1}&\text{if}\quad \alpha>0, \ p=1+\alpha/(n-\alpha),\\
		\varepsilon^{-(p-1)}&\text{if}\quad \alpha>0,\ 1<p<1+\alpha/(n-\alpha)
	\end{cases}
$$
for any $\varepsilon\in (0,1]$,
where
$$
	\kappa=\frac{2(1+\beta)}{2-\alpha}\left(\frac{1}{p-1}-\frac{n-\alpha}{2}\right).
$$
\end{theorem}

\begin{remark}
The results of Theorem 2.3 and Proposition 2.2
can be expressed by the following table:

\newlength{\myheight}
\setlength{\myheight}{1cm}
\newlength{\myheighta}
\setlength{\myheighta}{3cm}
\begin{table}[h]
  \begin{tabular}{|c|c|c|} \hline
    \ &$\alpha=0$&
	$\beta=0$\\ \hline
    \parbox[c][\myheight][c]{0cm}{} $p_c$&
	$\displaystyle 1+\frac{2}{n}$&
	$\displaystyle 1+\frac{2}{n-\alpha}$
	\\ \hline
    \parbox[c][\myheighta][c]{0cm}{} $T_{\varepsilon}\lesssim $&
	$\varepsilon^{-1/\kappa}$&
    	$\left\{ \begin{array}{lll}
    		\varepsilon^{-1/\kappa},\quad
			\displaystyle \left(1+\frac{\alpha}{n-\alpha}<p<1+\frac{2}{n-\alpha}\right)\\[7pt]
		\varepsilon^{-(p-1)}(\log(\varepsilon^{-1}))^{p-1},\quad
			\displaystyle \left(p=1+\frac{\alpha}{n-\alpha}\right)\\[7pt]
		\varepsilon^{-(p-1)},\quad 
			\displaystyle \left(1<p<1+\frac{\alpha}{n-\alpha}\right)
		\end{array}\right.$
	\\ \hline
	\parbox[c][\myheight][c]{0cm}{}
		$T_{\varepsilon}\gtrsim$&
		$\varepsilon^{-1/\kappa+\delta}$&
		$\varepsilon^{-1/\kappa+\delta}$
	\\ \hline
    \parbox[c][\myheight][c]{0cm}{} $\kappa$&
	$\displaystyle (1+\beta)\left(\frac{1}{p-1}-\frac{n}{2}\right)$&
	$\displaystyle \frac{2}{2-\alpha}\left(\frac{1}{p-1}-\frac{n-\alpha}{2}\right)$
	\\ \hline
  \end{tabular}
\end{table}

\end{remark}

\begin{remark}
It is expected that the rate $\kappa$ in Theorems 2.3 is sharp
except for the case $\alpha>0, 1<p\le 1+\alpha/(n-\alpha)$
from Proposition 2.2.
\end{remark}

\begin{remark}
The explicit form of $\Phi=\langle x\rangle^{-\alpha}(1+t)^{-\beta}$ is not necessary.
Indeed, we can treat more general coefficients,
for example,
$\Phi(t,x)=a(x)$
satisfying $a\in C(\mathbf{R}^n)$ and $0\le a(x)\lesssim \langle x\rangle^{-\alpha}$,
or
$\Phi(t,x)=b(t)$
satisfying $b\in C^1([0,\infty))$ and $b(t)\sim (1+t)^{-\beta}$.
\end{remark}

\begin{remark}
The same conclusion of Theorem 2.3
is valid for the corresponding heat equation
$-\Delta v+\Phi(t,x)v_t=|v|^p$
in the same manner as our proof.
\end{remark}

Our proof is based on a test function method.
Zhang \cite{Z} also used a similar way to determine the critical exponent
for the case $\alpha=\beta=0$.
By using his method, many blow-up results were obtained
for variable coefficient cases
(see \cite{DL, ITY, LNZ}).
However, 
the method of \cite{Z} was based on a contradiction argument
and so upper estimates of the lifespan cannot be obtained.
To avoid the contradiction argument,
we use an idea by Kuiper \cite{K}.
He obtained an upper bound of the lifespan for some parabolic equations
(see also \cite{I, Sun}).
We note that to treat the time-dependent damping case,
we also use a transformation of equation by Lin, Nishihara and Zhai \cite{LNZ}
(see also \cite{DL}).

At the end of this section, we explain some notation and terminology
used throughout this paper.
We put
$$
	\|f\|_{L^p(\mathbf{R}^n)}:=
	\left(\int_{\mathbf{R}^n}|f(x)|^pdx\right)^{1/p}.
$$
We denote the usual Sobolev space by $H^1(\mathbf{R}^n)$.
For an interval $I$ and a Banach space $X$, we define $C^r(I;X)$ as the Banach space whose
element is an $r$-times continuously differentiable mapping from $I$ to $X$ with respect to the topology in $X$.
The letter $C$ indicates the generic constant, which may change from line to line. 
We also use the symbols $\lesssim$ and $\sim$.
The relation $f\lesssim g$ means $f\le Cg$ with some constant $C>0$
and $f\sim g$ means $f\lesssim g$ and $g\lesssim f$.


\section{Proof of Theorem 2.3}
We first note that
if $T_{\varepsilon}\le C$,
where $C$ is a positive constant
depending only on $n, p, \alpha, \beta, u_0, u_1$,
then it is obvious that
$T_{\varepsilon}\le C\varepsilon^{-1/\kappa}$
for any $\kappa>0$ and $\varepsilon\in (0,1]$.
Therefore, once a constant $C=C(n, p, \alpha, \beta, u_0, u_1)$ is given,
we may assume that $T_{\varepsilon}> C$.

In the case $\beta\neq 0$, (\ref{eq11}) is not divergence form
and so we cannot apply the test function method.
Therefore, we need to transform the equation (\ref{eq11}) into divergence form.
The following idea was introduced by Lin, Nishihara and Zhai \cite{LNZ}.
Let $g(t)$ be the solution of the ordinary differential equation
$$
	\begin{cases}
		-g^{\prime}(t)+(1+t)^{-\beta}g(t)=1,\\
		g(0)=B^{-1}.
	\end{cases}
$$
The solution $g(t)$ is explicitly given by
$$
	g(t)=e^{\int_0^t(1+s)^{-\beta}ds}\left(B^{-1}-\int_0^te^{-\int_0^{\tau}(1+s)^{-\beta}ds}d\tau\right).
$$
By the de l'H${\rm \hat{o}}$pital theorem, we have
$$
	\lim_{t\to\infty}(1+t)^{-\beta}g(t)=1
$$
and so $g(t)\sim (1+t)^{\beta}$.
We note that $B=1$ and $g(t)\equiv 1$ if $\beta=0$.
By the definition of $g(t)$, we also have $|g^{\prime}(t)|\lesssim |(1+t)^{-\beta}g(t)-1|\lesssim 1$.
Multiplying the equation (\ref{eq11}) by $g(t)$, we obtain the divergence form
\begin{equation}
\label{eq21}
	(gu)_{tt}-\Delta (gu)-((g^{\prime}-1)\langle x\rangle^{-\alpha}u)_t=g|u|^p,
\end{equation}
here we note that $\alpha\beta=0$.
Therefore, we can apply the test function method to (\ref{eq21}).

We introduce the following test functions:
\begin{align*}
	\phi(x)&:=\begin{cases}
		\exp(-1/(1-|x|^2))&(|x|<1),\\
		0&(|x|\ge 1),
	\end{cases}\\
	\eta(t)&:=\begin{cases}
		\displaystyle \frac{\exp(-1/(1-t^2))}{\exp(-1/(t^2-1/4))+\exp(-1/(1-t^2))}&(1/2<t<1),\\
		1&(0\le t\le 1/2),\\
		0&(t\ge 1).
	\end{cases}
\end{align*}
It is obvious that $\phi\in C_0^{\infty}(\mathbf{R}^n), \eta\in C_0^{\infty}([0,\infty))$
and there exists a constant $C>0$ such that for all $|x|<1$ we have
$$
	\frac{|\nabla \phi(x)|^2}{\phi(x)}\le C.
$$
Using this estimate, we can prove that there exists a constant $C>0$ such that
the estimate
\begin{equation}
\label{eq22}
	|\Delta \phi(x)|\le C\phi(x)^{1/p}
\end{equation}
is true for all $|x|<1$.
Indeed, putting $\varphi:=\phi^{1/q}$ with $q=p/(p-1)$, we have
$$
	|\Delta \phi(x)|=|\Delta (\varphi(x)^q)|
	\lesssim |\Delta \varphi(x)|\varphi(x)^q+|\nabla \varphi(x)|^2\varphi(x)^{q-2}
	\lesssim \varphi(x)^{q-1}=\phi(x)^{1/p}.
$$
In the same way, we can also prove that
\begin{equation}
\label{eq23}
	|\eta^{\prime}(t)|\le C\eta(t)^{1/p},\quad |\eta^{\prime\prime}(t)|\le C\eta(t)^{1/p}
\end{equation}
for $t\in [0,1)$.

Let
$u$
be a solution on $[0,T_{\varepsilon})$
and
$\tau\in (\tau_0, T_{\varepsilon}), R\ge R_0$
parameters, where $\tau_0\in [1,T_{\varepsilon}), R_0>0$ are defined later.
We define
$$
	\psi_{\tau, R}(t,x):=\eta_{\tau}(t)\phi_R(x):=\eta(t/\tau)\phi(x/R)
$$
and
\begin{align*}
	I_{\tau,R}&:=\int_{[0,\tau)\times B_R}g(t)|u(t,x)|^p\psi_{\tau,R}(t,x)dxdt,\\
	J_R&:=\varepsilon\int_{B_R}(\langle x\rangle^{-\alpha}Bu_0(x)+u_1(x))\phi_R(x)dx,
\end{align*}
where $B_R=\{ |x|< R \}$.
Since $\psi_{\tau,R}\in C_0^{\infty}([0,T_{\varepsilon})\times \mathbf{R}^n)$
and $u$ is a solution on $[0,T_{\varepsilon})$,
we have
\begin{align*}
	I_{\tau,R}+J_R&=\int_{[0,\tau)\times B_R}g(t)u\partial_t^2\psi_{\tau,R}dxdt
		-\int_{[0,\tau)\times B_R}g(t)u\Delta\psi_{\tau,R}dxdt\\
		&\quad+\int_{[0,\tau)\times B_R}(g^{\prime}(t)-1)\langle x\rangle^{-\alpha}u\partial_t\psi_{\tau,R}dxdt\\
		&=:K_1+K_2+K_3.
\end{align*}
Here we have used the property $\partial_t\psi(0,x)=0$
and substituted the test function $g(t)\psi(t,x)$ into the definition of solution (\ref{sol}).
We note that for the corresponding heat equation,
we have the same decomposition without the term $K_1$
and so we can obtain the same conclusion
(see Remark 2.4).
We first estimate $K_1$.
By the H\"older inequality and (\ref{eq23}), we have
\begin{align}
\label{eq24}
	K_1&\le
		\tau^{-2}\int_{[0,\tau)\times B_R}g(t)|u||\eta^{\prime\prime}(t/\tau)|\phi_R(x)dxdt\\
	&\le C\tau^{-2}\int_{[\tau/2,\tau)\times B_R}g(t)|u|\eta(t)^{1/p}\phi_R(x)dxdt\notag\\
	&\le\tau^{-2}I_{\tau,R}^{1/p}
		\left(\int_{\tau/2}^{\tau}g(t)dt
		\cdot\int_{B_R}\phi_R(x)dx\right)^{1/q}\notag\\
	&\le C\tau^{-2+1/q}(1+\tau)^{\beta/q}R^{n/q}I_{\tau,R}^{1/p}.\notag
\end{align}
Using (\ref{eq22}) and a similar calculation, we obtain
\begin{align}
\label{eq25}
	K_2&\le R^{-2}\int_{[0,\tau)\times B_R}g(t)|u||\Delta\phi(x/R)|\eta(t/\tau)dxdt\\
	&\le CR^{-2}\int_{[0,\tau)\times B_R}g(t)|u||\phi(x/R)|^{1/p}\eta(t/\tau)dxdt\notag\\
	&\le CR^{-2}I_{\tau,R}^{1/p}
		\left(\int_0^{\tau}g(t)\eta(t/\tau)dt\cdot\int_{B_R}1dx\right)^{1/q}\notag\\
	&\le C(1+\tau)^{(1+\beta)/q}R^{-2+n/q}I_{\tau,R}^{1/p}.\notag
\end{align}
For $K_3$, using (\ref{eq23}) and $|g^{\prime}(t)-1|\lesssim C$, we have
\begin{align}
\label{eq26}
	K_3&\le \tau^{-1}\int_{[0,\tau)\times B_R}\langle x\rangle^{-\alpha}|u||\eta^{\prime}(t/\tau)|\phi_R(x)dxdt\\
	&\le \tau^{-1}I_{\tau,R}^{1/p}
		\left(\int_{\tau/2}^{\tau}g(t)^{-q/p}dt
		\cdot\int_{B_R}\langle x\rangle^{-\alpha q}\phi_R(x)dx\right)^{1/q}\notag\\
	&\le C\tau^{-1+1/q}(1+\tau)^{-\beta/p}F_{p,\alpha}(R)I_{\tau,R}^{1/p},\notag
\end{align}
where
$$
	F_{p,\alpha}(R)=\begin{cases}
		R^{-\alpha+n/q}&(\alpha q<n),\\
		(\log(1+R))^{1/q}&(\alpha q=n),\\
		1&(\alpha q>n).
	\end{cases}
$$
Thus, putting
$$
	D(\tau,R):=\tau^{-(1+\beta)/p}(\tau^{-1+\beta}R^{q/n}+\tau^{1+\beta}R^{-2+q/n}+F_{p,\alpha}(R))
$$
and combining this with the estimates (\ref{eq24})-(\ref{eq26}),
we have
\begin{equation}
\label{eq27}
	J_R\le CD(\tau,R)I_{\tau,R}^{1/p}-I_{\tau,R}.
\end{equation}
Now we use a fact that the inequality
$$
	ac^b-c\le (1-b)b^{b/(1-b)}a^{1/(1-b)}
$$
holds for all $a>0, 0<b<1, c\ge 0$.
We can immediately prove it by considering
the maximal value of the function $f(c)=ac^b-c$.
From this and (\ref{eq27}), we obtain
\begin{equation}
\label{eq28}
	J_R\le CD(\tau, R)^q.
\end{equation}
On the other hand, by the assumption on the data and
the Lebesgue dominated convergence theorem,
there exist $C>0$ and $R_0$ such that
$J_R\ge C\varepsilon$ holds for all $R>R_0$.
Combining this with (\ref{eq28}), we have
\begin{equation}
\label{eq29}
	\varepsilon\le CD(\tau, R)^q
\end{equation}
for all $\tau\in (\tau_0,T_{\varepsilon})$ and $R>R_0$.
Now we difine
$$
	\tau_0:=\max\{ 1, R_0^{(2-\alpha)/(1+\beta)} \},
$$
and we substitute
\begin{equation}
\label{eq210}
	R=\begin{cases}
		\tau^{(1+\beta)/(2-\alpha)}&(\alpha q<n),\\
		\tau&(\alpha q\ge n)
	\end{cases}
\end{equation}
into (\ref{eq29}).
Here we note that $R>R_0$ if $R$ is given by (\ref{eq210}).
As was mentioned at the beginning of this section,
we may assume that
$\tau_0<T_{\varepsilon}$.
Finally, we have
$$
	\varepsilon\le C\begin{cases}
		\tau^{-\kappa}&(\alpha q<n),\\
		\tau^{-1/(p-1)}\log(1+\tau)&(\alpha q=n),\\
		\tau^{-1/(p-1)}&(\alpha q>n),
	\end{cases}
$$
with
$$
	\kappa=\frac{2(1+\beta)}{2-\alpha}\left(\frac{1}{p-1}-\frac{n-\alpha}{2}\right).
$$
We can rewrite this relation as
$$
	\tau \le C\begin{cases}
		\varepsilon^{-1/\kappa}&\text{if}\quad 1+\alpha/(n-\alpha)<p<1+2/(n-\alpha),\\
		\varepsilon^{-(p-1)}(\log(\varepsilon^{-1}))^{p-1}&\text{if}\quad \alpha>0, \ p=1+\alpha/(n-\alpha),\\
		\varepsilon^{-(p-1)}&\text{if}\quad \alpha>0,\ 1<p<1+\alpha/(n-\alpha).
	\end{cases}
$$
Here we note that $\kappa>0$ if and only if $1<p<1+2/(n-\alpha)$
and that $\alpha q=n$ is equivalent to $p=1+\alpha/(n-\alpha)$.
Since $\tau$ is arbitrary in $(\tau_0, T_{\varepsilon})$,
we can obtain the conclusion of the theorem.

\section*{Acknowledgement}
The authors are deeply grateful to Professor Tatsuo Nishitani.
He gave us helpful advice.
They would also like to express their deep gratitude to
an anonymous referee for many useful suggestions and comments.

\end{document}